\documentclass{amsart}
\usepackage{amsfonts}

\newtheorem{theorem}{Theorem}
\theoremstyle{plain}

\newtheorem{corollary}{Corollary}

\newtheorem{definition}{Definition}
\newtheorem{example}{Example}

\numberwithin{equation}{section}
\input{tcilatex}

\begin{document}
\title{Smarandache Curves According to Bishop Frame in Euclidean 3-Space}
\author{Muhammed \c{C}etin}
\address{Usak University, Instutition of Science and Technology, MSc
Student, 1 Eylul Campus, 64200, Usak-TURKEY}
\email{mat.mcetin@hotmail.com}
\author{Y\i lmaz Tun\c{c}er}
\address{Usak University, Faculty of Sciences and Arts, Department of
Mathematics, 1 Eylul Campus, 64200, Usak-TURKEY}
\email{yilmaz.tuncer@usak.edu.tr}
\author{Murat Kemal Karacan}
\address{Usak University, Faculty of Sciences and Arts, Department of
Mathematics, 1 Eylul Campus, 64200 Usak-TURKEY}
\email{murat.karacan@usak.edu.tr}

\begin{abstract}
In this paper, we investigate special Smarandache curves according to Bishop
frame in Euclidean 3-space and we give some differential geometric
properties of Smarandache curves. Also we find the centers of the osculating
spheres and curvature spheres of Smarandache curves.
\end{abstract}

\maketitle

\section{\protect\bigskip Introduction}

Special Smarandache curves have been investigated by some differential
geometers [1,8]. A regular curve in Minkowski space-time, whose position
vector is composed by Frenet frame vectors on another regular curve, is
called a \textit{Smarandache\ Curve}. M. Turgut and S. Yilmaz have defined a
special case of such curves and call it \textit{Smarandache }$TB_{2}$ 
\textit{Curves} in the space $E_{1}^{4}$ [8]. They have dealed with a
special Smarandache curves which is defined by the tangent and second
binormal vector fields. They have called such curves as \textit{Smarandache }%
$TB_{2}$\textit{\ Curves}. Additionally, they have computed formulas of this
kind curves by the method expressed in [10]. A. T. Ali has introduced some
special Smarandache curves in the Euclidean space. He has studied
Frenet-Serret invariants of a special case [1].

In this paper, we investigate special Smarandache curves such as \textit{%
Smarandache Curves }$TN_{1}$\textit{, }$TN_{2}$\textit{, }$N_{1}N_{2}$%
\textit{\ and }$TN_{1}N_{2}$ according to Bishop frame in Euclidean 3-space.
Furthermore, we find differential geometric properties of these special
curves and we calculate first and second curvature (natural curvatures) of
these curves. Also we find the centers of the curvature spheres and
osculating spheres of Smarandache curves.

\section{Preliminaries}

\bigskip The Bishop frame or parallel transport frame is an alternative
approach to defining a moving frame that is well defined even when the curve
has vanishing second derivative. We can parallel transport an orthonormal
frame along a curve simply by parallel transporting each component of the
frame. The parallel transport frame is based on the observation that, while $%
T(s)$ for a given curve model is unique, we may choose any convenient
arbitrary basis $\left( N_{1}(s),N_{2}(s)\right) $ for the remainder of the
frame, so long as it is in the normal plane perpendicular to $T(s)$ at each
point. If the derivatives of $\left( N_{1}(s),N_{2}(s)\right) $ depend only
on $T(s)$ and not each other we can make $N_{1}(s)$ and $N_{2}(s)$ vary
smoothly throughout the path regardless of the curvature [2,6,7].

In addition, suppose the curve $\alpha $ is an arclength-parametrized $C^{2}$
curve. Suppose we have $C^{1}$ unit vector fields $N_{1}$ and $N_{2}=T\wedge 
$ $N_{1}$ along the curve $\alpha $ so that%
\begin{equation*}
\left\langle T,N_{1}\right\rangle =\left\langle T,N_{2}\right\rangle
=\left\langle N_{1},N_{2}\right\rangle =0,
\end{equation*}%
i.e., $T$ $,N_{1},N_{2}$ will be a smoothly varying right-handed orthonormal
frame as we move along the curve. (To this point, the Frenet frame would
work just fine if the curve were $C^{3}$ with $\kappa \neq 0$) But now we
want to impose the extra condition that $\left\langle N_{1}^{\prime
},N_{2}\right\rangle =0$. We say the unit first normal vector field $N_{1}$
is parallel along the curve $\alpha .$ This means that the only change of $%
N_{1}$ is in the direction of $T$. A Bishop frame can be defined even when a
Frenet frame cannot (e.g., when there are points with $\kappa =0$).
Therefore, we have the alternative frame equations%
\begin{equation*}
\left[ 
\begin{array}{l}
T\,^{\prime } \\ 
N\,_{1}^{\prime } \\ 
N\,_{2}^{\prime }%
\end{array}%
\right] =\left[ 
\begin{array}{ccc}
0 & k_{1} & k_{2} \\ 
-k_{1} & 0 & 0 \\ 
-k_{2} & 0 & 0%
\end{array}%
\right] \left[ 
\begin{array}{l}
T \\ 
N_{1} \\ 
N_{2}%
\end{array}%
\right] \text{.}
\end{equation*}

One can show that%
\begin{equation*}
\kappa (s)=\sqrt{k_{1}^{2}+k_{2}^{2}},\text{ }\theta (s)=\arctan \left( 
\frac{k_{2}}{k_{1}}\right) ,\text{ }k_{1}\neq 0,~\tau (s)=-\frac{d\theta (s)%
}{ds}
\end{equation*}
so that $k_{1}$ and $k_{2}$ effectively correspond to a Cartesian coordinate
system for the polar coordinates $\kappa ,$ $\theta $ with $\theta =-\int
\tau (s)ds$. The orientation of the parallel transport frame includes the
arbitrary choice of integration constant $\theta _{0}$, which disappears
from $\tau $ (and hence from the Frenet frame) due to the differentiation
[2,6,7].

Bishop curvatures are defined by%
\begin{equation}
k_{1}=\kappa \cos \theta \text{, \ \ \ }k_{2}=\kappa \sin \theta  \label{x3}
\end{equation}%
and%
\begin{equation}
T=T\text{ \ },\text{ \ \ \ \ }N_{1}=N\cos \theta -B\sin \theta \text{ \ },%
\text{ \ \ \ \ }N_{2}=N\sin \theta +B\cos \theta  \label{x4}
\end{equation}%
[9].

\section{Smarandache Curves According to Bishop Frame}

In [1], author gave following definitions:

\begin{definition}
Let $\alpha =\alpha (s)$ be a unit speed regular curve in $E^{3}$ and $%
\left\{ T,N,B\right\} $ be its moving Serret-Frenet frame. Smarandache
curves $TN$ are defined by%
\begin{equation*}
\beta (s^{\ast })=\frac{1}{\sqrt{2}}\left( T+N\right) .
\end{equation*}
\end{definition}

\begin{definition}
Let $\alpha =\alpha (s)$ be a unit speed regular curve in $E^{3}$ and $%
\left
\{ T,N,B\right \} $ be its moving Serret-Frenet frame. Smarandache
curves $NB $ are defined by%
\begin{equation*}
\beta (s^{\ast })=\frac{1}{\sqrt{2}}\left( N+B\right) .
\end{equation*}
\end{definition}

\begin{definition}
Let $\alpha =\alpha (s)$ be a unit speed regular curve in $E^{3}$ and $%
\left
\{ T,N,B\right \} $ be its moving Serret-Frenet frame. Smarandache
curves $TNB$ are defined by%
\begin{equation*}
\beta (s^{\ast })=\frac{1}{\sqrt{3}}\left( T+N+B\right) .
\end{equation*}
\end{definition}

In this section, we investigate Smarandache curves according to Bishop frame
in Euclidean 3-space.\ Let $\alpha =\alpha (s)$ be a unit speed regular
curve in $E^{3}$ and denote by $\left\{ T_{\alpha },N_{1}^{\alpha
},N_{2}^{\alpha }\right\} $ the moving Bishop frame along the curve $\alpha $%
. The following Bishop formulae is given by%
\begin{equation}
\overset{\cdot }{T_{\alpha }}=k_{1}^{\alpha }N_{1}^{\alpha }+k_{2}^{\alpha
}N_{2}^{\alpha }\text{ \ , \ }\overset{\cdot }{N_{1}^{\alpha }}%
=-k_{1}^{\alpha }T_{\alpha }\text{ \ , \ }\overset{\cdot }{N_{2}^{\alpha }}%
=-k_{2}^{\alpha }T_{\alpha }  \label{x}
\end{equation}%
with%
\begin{equation*}
T_{\alpha }=N_{1}^{\alpha }\wedge N_{2}^{\alpha }\ ,\ N_{1}^{\alpha
}=-T_{\alpha }\wedge N_{2}^{\alpha }\ ,\ N_{2}^{\alpha }=T_{\alpha }\wedge
N_{1}^{\alpha }.
\end{equation*}

\subsection{$TN_{1}-$Smarandache Curves}

\begin{definition}
Let $\alpha =\alpha (s)$ be a unit speed regular curve in $E^{3}$ and $%
\left\{ T_{\alpha },N_{1}^{\alpha },N_{2}^{\alpha }\right\} $ be its moving
Bishop frame. $T_{\alpha }N_{1}^{\alpha }-$Smarandache curves can be defined
by%
\begin{equation}
\beta (s^{\ast })=\frac{1}{\sqrt{2}}\left( T_{\alpha }+N_{1}^{\alpha
}\right) .  \label{a1}
\end{equation}
\end{definition}

Now, we can investigate Bishop invariants of $T_{\alpha }N_{1}^{\alpha }-$%
Smarandache curves according to $\alpha =\alpha (s).$ Differentiating (\ref%
{a1}) with respect to $s,$ we get%
\begin{equation}
\overset{\cdot }{\beta }=\frac{d\beta }{ds^{\ast }}\frac{ds^{\ast }}{ds}=%
\frac{-1}{\sqrt{2}}\left( k_{1}^{\alpha }T_{\alpha }-k_{1}^{\alpha
}N_{1}^{\alpha }-k_{2}^{\alpha }N_{2}^{\alpha }\right)  \label{a5}
\end{equation}%
and%
\begin{equation*}
T_{\beta }\frac{ds^{\ast }}{ds}=\frac{-1}{\sqrt{2}}\left( k_{1}^{\alpha
}T_{\alpha }-k_{1}^{\alpha }N_{1}^{\alpha }-k_{2}^{\alpha }N_{2}^{\alpha
}\right)
\end{equation*}%
where%
\begin{equation}
\frac{ds^{\ast }}{ds}=\sqrt{\frac{2\left( k_{1}^{\alpha }\right) ^{2}+\left(
k_{2}^{\alpha }\right) ^{2}}{2}}.  \label{a3}
\end{equation}%
The tangent vector of\ curve $\beta $ can be written as follow,%
\begin{equation}
T_{\beta }=\frac{-1}{\sqrt{2\left( k_{1}^{\alpha }\right) ^{2}+\left(
k_{2}^{\alpha }\right) ^{2}}}\left( k_{1}^{\alpha }T_{\alpha }-k_{1}^{\alpha
}N_{1}^{\alpha }-k_{2}^{\alpha }N_{2}^{\alpha }\right) .  \label{a2}
\end{equation}%
Differentiating (\ref{a2}) with respect to $s$, we obtain%
\begin{equation}
\frac{dT_{\beta }}{ds^{\ast }}\frac{ds^{\ast }}{ds}=\frac{1}{\left( 2\left(
k_{1}^{\alpha }\right) ^{2}+\left( k_{2}^{\alpha }\right) ^{2}\right) ^{%
\frac{3}{2}}}\left( \lambda _{1}T_{\alpha }+\lambda _{2}N_{1}^{\alpha
}+\lambda _{3}N_{2}^{\alpha }\right)  \label{a4}
\end{equation}%
where%
\begin{eqnarray*}
\lambda _{1} &=&-\overset{\cdot }{k_{1}^{\alpha }}\left( k_{2}^{\alpha
}\right) ^{2}-2\left( k_{1}^{\alpha }\right) ^{4}-3\left( k_{1}^{\alpha
}\right) ^{2}\left( k_{2}^{\alpha }\right) ^{2}-\left( k_{2}^{\alpha
}\right) ^{4}+k_{1}^{\alpha }k_{2}^{\alpha }\overset{\cdot }{k_{2}^{\alpha }}
\\
\lambda _{2} &=&-2\left( k_{1}^{\alpha }\right) ^{4}-\left( k_{1}^{\alpha
}\right) ^{2}\left( k_{2}^{\alpha }\right) ^{2}+\overset{\cdot }{%
k_{1}^{\alpha }}\left( k_{2}^{\alpha }\right) ^{2}-k_{1}^{\alpha
}k_{2}^{\alpha }\overset{\cdot }{k_{2}^{\alpha }} \\
\lambda _{3} &=&-2\left( k_{1}^{\alpha }\right) ^{3}k_{2}^{\alpha
}-k_{1}^{\alpha }\left( k_{2}^{\alpha }\right) ^{3}+2\left( k_{1}^{\alpha
}\right) ^{2}\overset{\cdot }{k_{2}^{\alpha }}-2k_{1}^{\alpha }\overset{%
\cdot }{k_{1}^{\alpha }}k_{2}^{\alpha }.
\end{eqnarray*}%
Substituting (\ref{a3}) in (\ref{a4}), we get%
\begin{equation*}
T_{\beta }^{\prime }=\frac{\sqrt{2}}{\left( 2\left( k_{1}^{\alpha }\right)
^{2}+\left( k_{2}^{\alpha }\right) ^{2}\right) ^{2}}\left( \lambda
_{1}T_{\alpha }+\lambda _{2}N_{1}^{\alpha }+\lambda _{3}N_{2}^{\alpha
}\right) .
\end{equation*}%
Then, the curvature and principal normal vector field of curve $\beta $ are
respectively,%
\begin{equation*}
\kappa _{\beta }=\left\Vert T_{\beta }^{\prime }\right\Vert =\frac{\sqrt{2}%
\sqrt{\lambda _{1}^{2}+\lambda _{2}^{2}+\lambda _{3}^{2}}}{\left( 2\left(
k_{1}^{\alpha }\right) ^{2}+\left( k_{2}^{\alpha }\right) ^{2}\right) ^{2}}
\end{equation*}%
and%
\begin{equation*}
N_{\beta }=\frac{1}{\sqrt{\lambda _{1}^{2}+\lambda _{2}^{2}+\lambda _{3}^{2}}%
}\left( \lambda _{1}T_{\alpha }+\lambda _{2}N_{1}^{\alpha }+\lambda
_{3}N_{2}^{\alpha }\right) .
\end{equation*}%
On the other hand, we express%
\begin{equation*}
B_{\beta }=\frac{1}{\sqrt{2\left( k_{1}^{\alpha }\right) ^{2}+\left(
k_{2}^{\alpha }\right) ^{2}}\sqrt{\lambda _{1}^{2}+\lambda _{2}^{2}+\lambda
_{3}^{2}}}\left\vert 
\begin{array}{ccc}
T_{\alpha } & N_{1}^{\alpha } & N_{2}^{\alpha } \\ 
-k_{1}^{\alpha } & k_{1}^{\alpha } & k_{2}^{\alpha } \\ 
\lambda _{1} & \lambda _{2} & \lambda _{3}%
\end{array}%
\right\vert .
\end{equation*}%
So, the binormal vector of curve $\beta $ is%
\begin{equation*}
B_{\beta }=\frac{1}{\sqrt{2\left( k_{1}^{\alpha }\right) ^{2}+\left(
k_{2}^{\alpha }\right) ^{2}}\sqrt{\lambda _{1}^{2}+\lambda _{2}^{2}+\lambda
_{3}^{2}}}\left( \sigma _{1}T_{\alpha }+\sigma _{2}N_{1}^{\alpha }+\sigma
_{3}N_{2}^{\alpha }\right)
\end{equation*}%
where%
\begin{equation*}
\sigma _{1}=\lambda _{3}k_{1}^{\alpha }-\lambda _{2}k_{2}^{\alpha }\text{ \
, \ }\sigma _{2}=\lambda _{3}k_{1}^{\alpha }+\lambda _{1}k_{2}^{\alpha }%
\text{ \ , \ }\sigma _{3}=-k_{1}^{\alpha }\left( \lambda _{1}+\lambda
_{2}\right) .
\end{equation*}%
We differentiate (\ref{a5}) with respect to $s$ in order to calculate the
torsion%
\begin{equation*}
\overset{\cdot \cdot }{\beta }=\frac{-1}{\sqrt{2}}\left\{ (\overset{\cdot }{%
k_{1}^{\alpha }}+\left( k_{1}^{\alpha }\right) ^{2}+\left( k_{2}^{\alpha
}\right) ^{2})T_{\alpha }-(\overset{\cdot }{k_{1}^{\alpha }}-\left(
k_{1}^{\alpha }\right) ^{2})N_{1}^{\alpha }-(\overset{\cdot }{k_{2}^{\alpha }%
}-k_{1}^{\alpha }k_{2}^{\alpha })N_{2}^{\alpha }\right\} ,
\end{equation*}%
and similarly 
\begin{equation*}
\overset{\cdot \cdot \cdot }{\beta }=\frac{1}{\sqrt{2}}\left( \rho
_{1}T_{\alpha }+\rho _{2}N_{1}^{\alpha }+\rho _{3}N_{2}^{\alpha }\right) ,
\end{equation*}%
where 
\begin{eqnarray*}
\rho _{1} &=&-\overset{\cdot \cdot }{k_{1}^{\alpha }}-3k_{1}^{\alpha }%
\overset{\cdot }{k_{1}^{\alpha }}-3k_{2}^{\alpha }\overset{\cdot }{%
k_{2}^{\alpha }}+\left( k_{1}^{\alpha }\right) ^{3}+k_{1}^{\alpha }\left(
k_{2}^{\alpha }\right) ^{2} \\
\rho _{2} &=&-3k_{1}^{\alpha }\overset{\cdot }{k_{1}^{\alpha }}-\left(
k_{1}^{\alpha }\right) ^{3}-k_{1}^{\alpha }\left( k_{2}^{\alpha }\right)
^{2}+\overset{\cdot \cdot }{k_{1}^{\alpha }} \\
\rho _{3} &=&-2\overset{\cdot }{k_{1}^{\alpha }}k_{2}^{\alpha }-\left(
k_{1}^{\alpha }\right) ^{2}k_{2}^{\alpha }-\left( k_{2}^{\alpha }\right)
^{3}-k_{1}^{\alpha }\overset{\cdot }{k_{2}^{\alpha }}+\overset{\cdot \cdot }{%
k_{2}^{\alpha }}.
\end{eqnarray*}%
The torsion of curve $\beta $ is

\begin{equation*}
\tau _{\beta }=\frac{\sqrt{2}\left \{ 
\begin{array}{c}
(\left( k_{1}^{\alpha }\right) ^{2}-\overset{\cdot }{k_{1}^{\alpha }})(\rho
_{3}k_{1}^{\alpha }+\rho _{1}k_{2}^{\alpha })+k_{1}^{\alpha }(\overset{\cdot 
}{k_{2}^{\alpha }}-k_{1}^{\alpha }k_{2}^{\alpha })(\rho _{1}+\rho _{2}) \\ 
-(\overset{\cdot }{k_{1}^{\alpha }}+\left( k_{1}^{\alpha }\right)
^{2}+\left( k_{2}^{\alpha }\right) ^{2})(\rho _{2}k_{2}^{\alpha }-\rho
_{3}k_{1}^{\alpha })%
\end{array}%
\right \} }{(k_{1}^{\alpha }\overset{\cdot }{k_{2}^{\alpha }}-\overset{\cdot 
}{k_{1}^{\alpha }}k_{2}^{\alpha })^{2}+(k_{1}^{\alpha }\overset{\cdot }{%
k_{2}^{\alpha }}-2\left( k_{1}^{\alpha }\right) ^{2}k_{2}^{\alpha }-\overset{%
\cdot }{k_{1}^{\alpha }}k_{2}^{\alpha }-\left( k_{2}^{\alpha }\right)
^{3})^{2}+(2\left( k_{1}^{\alpha }\right) ^{3}+k_{1}^{\alpha }\left(
k_{2}^{\alpha }\right) ^{2})^{2}}.
\end{equation*}

The first and second normal field of curve $\beta $ are as follow. Then,
from (\ref{x4}) we obtain

\begin{equation*}
N_{1}^{\beta }=\frac{1}{\sqrt{\mu }\eta }\left \{ 
\begin{array}{c}
\left( \sqrt{\mu }\cos \theta _{\beta }\lambda _{1}-\sin \theta _{\beta
}\sigma _{1}\right) T_{\alpha } \\ 
+\left( \sqrt{\mu }\cos \theta _{\beta }\lambda _{2}-\sin \theta _{\beta
}\sigma _{2}\right) N_{1}^{\alpha } \\ 
+\left( \sqrt{\mu }\cos \theta _{\beta }\lambda _{3}-\sin \theta _{\beta
}\sigma _{3}\right) N_{2}^{\alpha }%
\end{array}%
\right \} ,
\end{equation*}%
and

\begin{equation*}
N_{2}^{\beta }=\frac{1}{\sqrt{\mu }\eta }\left \{ 
\begin{array}{c}
\left( \sqrt{\mu }\sin \theta _{\beta }\lambda _{1}+\cos \theta _{\beta
}\sigma _{1}\right) T_{\alpha } \\ 
+\left( \sqrt{\mu }\sin \theta _{\beta }\lambda _{2}+\cos \theta _{\beta
}\sigma _{2}\right) N_{1}^{\alpha } \\ 
+\left( \sqrt{\mu }\sin \theta _{\beta }\lambda _{3}+\cos \theta _{\beta
}\sigma _{3}\right) N_{2}^{\alpha }%
\end{array}%
\right \}
\end{equation*}%
where $\theta _{\beta }=-\int \tau _{\beta }(s)ds$ and $\mu =2\left(
k_{1}^{\alpha }\right) ^{2}+\left( k_{2}^{\alpha }\right) ^{2}$, $\eta =%
\sqrt{\lambda _{1}^{2}+\lambda _{2}^{2}+\lambda _{3}^{2}}$

Now, we can calculate\ natural curvatures of curve $\beta $, so from (\ref%
{x3}) we get

\begin{equation*}
k_{1}^{\beta }=\frac{\sqrt{2(\lambda _{1}^{2}+\lambda _{2}^{2}+\lambda
_{3}^{2})}\cos \theta _{\beta }}{\left( 2\left( k_{1}^{\alpha }\right)
^{2}+\left( k_{2}^{\alpha }\right) ^{2}\right) ^{2}}
\end{equation*}%
and

\begin{equation*}
k_{2}^{\beta }=\frac{\sqrt{2(\lambda _{1}^{2}+\lambda _{2}^{2}+\lambda
_{3}^{2})}\sin \theta _{\beta }}{\left( 2\left( k_{1}^{\alpha }\right)
^{2}+\left( k_{2}^{\alpha }\right) ^{2}\right) ^{2}}.
\end{equation*}

\subsection{$TN_{2}-$Smarandache Curves}

\begin{definition}
Let $\alpha =\alpha (s)$ be a unit speed regular curve in $E^{3}$ and $%
\left\{ T_{\alpha },N_{1}^{\alpha },N_{2}^{\alpha }\right\} $ be its moving
Bishop frame. $T_{\alpha }N_{2}^{\alpha }-$Smarandache curves can be defined
by%
\begin{equation}
\beta (s^{\ast })=\frac{1}{\sqrt{2}}\left( T_{\alpha }+N_{2}^{\alpha
}\right) .  \label{b1}
\end{equation}
\end{definition}

Now, we can investigate Bishop invariants of $T_{\alpha }N_{2}^{\alpha }-$%
Smarandache curves according to $\alpha =\alpha (s).$ Differentiating (\ref%
{b1}) with respect to $s$, we get%
\begin{equation}
\overset{\cdot }{\beta }=\frac{d\beta }{ds^{\ast }}\frac{ds^{\ast }}{ds}=%
\frac{-1}{\sqrt{2}}\left( k_{2}^{\alpha }T_{\alpha }-k_{1}^{\alpha
}N_{1}^{\alpha }-k_{2}^{\alpha }N_{2}^{\alpha }\right)  \label{b5}
\end{equation}%
and%
\begin{equation*}
T_{\beta }\frac{ds^{\ast }}{ds}=\frac{-1}{\sqrt{2}}\left( k_{2}^{\alpha
}T_{\alpha }-k_{1}^{\alpha }N_{1}^{\alpha }-k_{2}^{\alpha }N_{2}^{\alpha
}\right)
\end{equation*}%
where%
\begin{equation}
\frac{ds^{\ast }}{ds}=\sqrt{\frac{\left( k_{1}^{\alpha }\right) ^{2}+2\left(
k_{2}^{\alpha }\right) ^{2}}{2}}.  \label{b3}
\end{equation}%
The tangent vector of curve $\beta $ can be written as follow,%
\begin{equation}
T_{\beta }=\frac{-1}{\sqrt{\left( k_{1}^{\alpha }\right) ^{2}+2\left(
k_{2}^{\alpha }\right) ^{2}}}\left( k_{2}^{\alpha }T_{\alpha }-k_{1}^{\alpha
}N_{1}^{\alpha }-k_{2}^{\alpha }N_{2}^{\alpha }\right) .  \label{b2}
\end{equation}%
Differentiating (\ref{b2}) with respect to $s$, we obtain%
\begin{equation}
\frac{dT_{\beta }}{ds^{\ast }}\frac{ds^{\ast }}{ds}=\frac{1}{\left( \left(
k_{1}^{\alpha }\right) ^{2}+2\left( k_{2}^{\alpha }\right) ^{2}\right) ^{%
\frac{3}{2}}}\left( \lambda _{1}T_{\alpha }+\lambda _{2}N_{1}^{\alpha
}+\lambda _{3}N_{2}^{\alpha }\right)  \label{b4}
\end{equation}%
where%
\begin{eqnarray*}
\lambda _{1} &=&-\left( k_{1}^{\alpha }\right) ^{2}\overset{\cdot }{%
k_{2}^{\alpha }}-\left( k_{1}^{\alpha }\right) ^{4}-3\left( k_{1}^{\alpha
}\right) ^{2}\left( k_{2}^{\alpha }\right) ^{2}-2\left( k_{2}^{\alpha
}\right) ^{4}+k_{1}^{\alpha }\overset{\cdot }{k_{1}^{\alpha }}k_{2}^{\alpha }
\\
\lambda _{2} &=&-\left( k_{1}^{\alpha }\right) ^{3}k_{2}^{\alpha
}-2k_{1}^{\alpha }\left( k_{2}^{\alpha }\right) ^{3}+2\overset{\cdot }{%
k_{1}^{\alpha }}\left( k_{2}^{\alpha }\right) ^{2}-2k_{1}^{\alpha
}k_{2}^{\alpha }\overset{\cdot }{k_{2}^{\alpha }} \\
\lambda _{3} &=&-\left( k_{1}^{\alpha }\right) ^{2}\left( k_{2}^{\alpha
}\right) ^{2}-2\left( k_{2}^{\alpha }\right) ^{4}+\left( k_{1}^{\alpha
}\right) ^{2}\overset{\cdot }{k_{2}^{\alpha }}-k_{1}^{\alpha }\overset{\cdot 
}{k_{1}^{\alpha }}k_{2}^{\alpha }.
\end{eqnarray*}%
Substituting (\ref{b3}) in (\ref{b4}), we get%
\begin{equation*}
T_{\beta }^{\prime }=\frac{\sqrt{2}}{\left( \left( k_{1}^{\alpha }\right)
^{2}+2\left( k_{2}^{\alpha }\right) ^{2}\right) ^{2}}\left( \lambda
_{1}T_{\alpha }+\lambda _{2}N_{1}^{\alpha }+\lambda _{3}N_{2}^{\alpha
}\right) .
\end{equation*}%
Then, the first curvature and the principal normal vector field of curve $%
\beta $ are respectively,%
\begin{equation*}
\kappa _{\beta }=\left\Vert T_{\beta }^{\prime }\right\Vert =\frac{\sqrt{2}%
\sqrt{\lambda _{1}^{2}+\lambda _{2}^{2}+\lambda _{3}^{2}}}{\left( \left(
k_{1}^{\alpha }\right) ^{2}+2\left( k_{2}^{\alpha }\right) ^{2}\right) ^{2}}
\end{equation*}%
and%
\begin{equation*}
N_{\beta }=\frac{1}{\sqrt{\lambda _{1}^{2}+\lambda _{2}^{2}+\lambda _{3}^{2}}%
}\left( \lambda _{1}T_{\alpha }+\lambda _{2}N_{1}^{\alpha }+\lambda
_{3}N_{2}^{\alpha }\right) .
\end{equation*}%
On the other hand, we express%
\begin{equation*}
B_{\beta }=\frac{1}{\sqrt{\left( k_{1}^{\alpha }\right) ^{2}+2\left(
k_{2}^{\alpha }\right) ^{2}}\sqrt{\lambda _{1}^{2}+\lambda _{2}^{2}+\lambda
_{3}^{2}}}\left\vert 
\begin{array}{ccc}
T_{\alpha } & N_{1}^{\alpha } & N_{2}^{\alpha } \\ 
-k_{2}^{\alpha } & k_{1}^{\alpha } & k_{2}^{\alpha } \\ 
\lambda _{1} & \lambda _{2} & \lambda _{3}%
\end{array}%
\right\vert .
\end{equation*}%
So, the binormal vector of curve $\beta $ is%
\begin{equation*}
B_{\beta }=\frac{1}{\sqrt{\left( k_{1}^{\alpha }\right) ^{2}+2\left(
k_{2}^{\alpha }\right) ^{2}}\sqrt{\lambda _{1}^{2}+\lambda _{2}^{2}+\lambda
_{3}^{2}}}\left( \sigma _{1}T_{\alpha }+\sigma _{2}N_{1}^{\alpha }+\sigma
_{3}N_{2}^{\alpha }\right)
\end{equation*}%
where%
\begin{equation*}
\sigma _{1}=\lambda _{3}k_{1}^{\alpha }-\lambda _{2}k_{2}^{\alpha }\text{ \
, \ }\sigma _{2}=\left( \lambda _{3}+\lambda _{1}\right) k_{2}^{\alpha }%
\text{ \ , \ }\sigma _{3}=-\left( \lambda _{2}k_{2}^{\alpha }+\lambda
_{1}k_{1}^{\alpha }\right) .
\end{equation*}%
We differentiate (\ref{b5}) with respect to $s$ in order to calculate the
torsion of curve $\beta $%
\begin{equation*}
\overset{\cdot \cdot }{\beta }=\frac{-1}{\sqrt{2}}\left\{ (\overset{\cdot }{%
k_{2}^{\alpha }}+\left( k_{1}^{\alpha }\right) ^{2}+\left( k_{2}^{\alpha
}\right) ^{2})T_{\alpha }-(\overset{\cdot }{k_{1}^{\alpha }}-k_{1}^{\alpha
}k_{2}^{\alpha })N_{1}^{\alpha }-(\overset{\cdot }{k_{2}^{\alpha }}-\left(
k_{2}^{\alpha }\right) ^{2})N_{2}^{\alpha }\right\} ,
\end{equation*}%
and similarly%
\begin{equation*}
\overset{\cdot \cdot \cdot }{\beta }=\frac{1}{\sqrt{2}}\left( \rho
_{1}T_{\alpha }+\rho _{2}N_{1}^{\alpha }+\rho _{3}N_{2}^{\alpha }\right) ,
\end{equation*}%
where 
\begin{eqnarray*}
\rho _{1} &=&-\overset{\cdot \cdot }{k_{2}^{\alpha }}-3k_{1}^{\alpha }%
\overset{\cdot }{k_{1}^{\alpha }}-3k_{2}^{\alpha }\overset{\cdot }{%
k_{2}^{\alpha }}+\left( k_{1}^{\alpha }\right) ^{2}k_{2}^{\alpha }+\left(
k_{2}^{\alpha }\right) ^{3} \\
\rho _{2} &=&-2k_{1}^{\alpha }\overset{\cdot }{k_{2}^{\alpha }}-\left(
k_{1}^{\alpha }\right) ^{3}-k_{1}^{\alpha }\left( k_{2}^{\alpha }\right)
^{2}-\overset{\cdot }{k_{1}^{\alpha }}k_{2}^{\alpha }+\overset{\cdot \cdot }{%
k_{1}^{\alpha }} \\
\rho _{3} &=&-3k_{2}^{\alpha }\overset{\cdot }{k_{2}^{\alpha }}-\left(
k_{1}^{\alpha }\right) ^{2}k_{2}^{\alpha }-\left( k_{2}^{\alpha }\right)
^{3}+\overset{\cdot \cdot }{k_{2}^{\alpha }}.
\end{eqnarray*}%
The torsion of curve\ $\beta $ is

\begin{equation*}
\tau _{\beta }=\frac{\sqrt{2}\left \{ 
\begin{array}{c}
k_{2}^{\alpha }(k_{1}^{\alpha }k_{2}^{\alpha }-\overset{\cdot }{%
k_{1}^{\alpha }})(\rho _{3}+\rho _{1})+(\overset{\cdot }{k_{2}^{\alpha }}%
-\left( k_{2}^{\alpha }\right) ^{2})(\rho _{2}k_{2}^{\alpha }+\rho
_{1}k_{1}^{\alpha }) \\ 
-(\overset{\cdot }{k_{2}^{\alpha }}+\left( k_{1}^{\alpha }\right)
^{2}+\left( k_{2}^{\alpha }\right) ^{2})(\rho _{2}k_{2}^{\alpha }-\rho
_{3}k_{1}^{\alpha })%
\end{array}%
\right \} }{(k_{1}^{\alpha }\overset{\cdot }{k_{2}^{\alpha }}-\overset{\cdot 
}{k_{1}^{\alpha }}k_{2}^{\alpha })^{2}+(-2\left( k_{2}^{\alpha }\right)
^{3}-\left( k_{1}^{\alpha }\right) ^{2}k_{2}^{\alpha })^{2}+(2k_{1}^{\alpha
}\left( k_{2}^{\alpha }\right) ^{2}-\overset{\cdot }{k_{1}^{\alpha }}%
k_{2}^{\alpha }+k_{1}^{\alpha }\overset{\cdot }{k_{2}^{\alpha }}+\left(
k_{1}^{\alpha }\right) ^{3})^{2}}.
\end{equation*}

The first and second normal field of curve $\beta $ are as follow. Then,
from (\ref{x4}) we obtain

\begin{equation*}
N_{1}^{\beta }=\frac{1}{\sqrt{\mu }\eta }\left \{ 
\begin{array}{c}
\left( \sqrt{\mu }\cos \theta _{\beta }\lambda _{1}-\sin \theta _{\beta
}\sigma _{1}\right) T_{\alpha } \\ 
+\left( \sqrt{\mu }\cos \theta _{\beta }\lambda _{2}-\sin \theta _{\beta
}\sigma _{2}\right) N_{1}^{\alpha } \\ 
+\left( \sqrt{\mu }\cos \theta _{\beta }\lambda _{3}-\sin \theta _{\beta
}\sigma _{3}\right) N_{2}^{\alpha }%
\end{array}%
\right \} ,
\end{equation*}%
and

\begin{equation*}
N_{2}^{\beta }=\frac{1}{\sqrt{\mu }\eta }\left \{ 
\begin{array}{c}
\left( \sqrt{\mu }\sin \theta _{\beta }\lambda _{1}+\cos \theta _{\beta
}\sigma _{1}\right) T_{\alpha } \\ 
+\left( \sqrt{\mu }\sin \theta _{\beta }\lambda _{2}+\cos \theta _{\beta
}\sigma _{2}\right) N_{1}^{\alpha } \\ 
+\left( \sqrt{\mu }\sin \theta _{\beta }\lambda _{3}+\cos \theta _{\beta
}\sigma _{3}\right) N_{2}^{\alpha }%
\end{array}%
\right \} ,
\end{equation*}%
where $\theta _{\beta }=-\int \tau _{\beta }(s)ds$ and $\mu =\left(
k_{1}^{\alpha }\right) ^{2}+2\left( k_{2}^{\alpha }\right) ^{2}$, \ $\eta =%
\sqrt{\lambda _{1}^{2}+\lambda _{2}^{2}+\lambda _{3}^{2}}.$

Now, we can calculate\ natural curvatures of curve $\beta $, so from (\ref%
{x3}) we get

\begin{equation*}
k_{1}^{\beta }=\frac{\sqrt{2(\lambda _{1}^{2}+\lambda _{2}^{2}+\lambda
_{3}^{2})}\cos \theta _{\beta }}{\left( \left( k_{1}^{\alpha }\right)
^{2}+2\left( k_{2}^{\alpha }\right) ^{2}\right) ^{2}}
\end{equation*}%
and

\begin{equation*}
k_{2}^{\beta }=\frac{\sqrt{2(\lambda _{1}^{2}+\lambda _{2}^{2}+\lambda
_{3}^{2})}\sin \theta _{\beta }}{\left( \left( k_{1}^{\alpha }\right)
^{2}+2\left( k_{2}^{\alpha }\right) ^{2}\right) ^{2}}.
\end{equation*}

\subsection{$N_{1}N_{2}-$Smarandache Curves}

\begin{definition}
Let $\alpha =\alpha (s)$ be a unit speed regular curve in $E^{3}$ and $%
\left\{ T_{\alpha },N_{1}^{\alpha },N_{2}^{\alpha }\right\} $ be its moving
Bishop frame. $N_{1}^{\alpha }N_{2}^{\alpha }-$Smarandache curves can be
defined by%
\begin{equation}
\beta (s^{\ast })=\frac{1}{\sqrt{2}}\left( N_{1}^{\alpha }+N_{2}^{\alpha
}\right) .  \label{c1}
\end{equation}
\end{definition}

Now, we can investigate Bishop invariants of $N_{1}^{\alpha }N_{2}^{\alpha
}- $Smarandache curves according to $\alpha =\alpha (s).$ Differentiating (%
\ref{c1}) with respect to $s$, we get%
\begin{equation}
\overset{\cdot }{\beta }=\frac{d\beta }{ds^{\ast }}\frac{ds^{\ast }}{ds}=%
\frac{-\left( k_{1}^{\alpha }+k_{2}^{\alpha }\right) }{\sqrt{2}}T_{\alpha }
\label{c5}
\end{equation}%
and%
\begin{equation*}
T_{\beta }\frac{ds^{\ast }}{ds}=\frac{-\left( k_{1}^{\alpha }+k_{2}^{\alpha
}\right) }{\sqrt{2}}T_{\alpha }
\end{equation*}%
where%
\begin{equation}
\frac{ds^{\ast }}{ds}=\frac{k_{1}^{\alpha }+k_{2}^{\alpha }}{\sqrt{2}}.
\label{c3}
\end{equation}%
The tangent vector of curve $\beta $ can be written as follow,%
\begin{equation}
T_{\beta }=-T_{\alpha }.  \label{c2}
\end{equation}%
Differentiating (\ref{c2}) with respect to $s$, we obtain%
\begin{equation}
\frac{dT_{\beta }}{ds^{\ast }}\frac{ds^{\ast }}{ds}=-k_{1}^{\alpha
}N_{1}^{\alpha }-k_{2}^{\alpha }N_{2}^{\alpha }  \label{c4}
\end{equation}%
Substituting (\ref{c3}) in (\ref{c4}), we get%
\begin{equation*}
T_{\beta }^{\prime }=\frac{-\sqrt{2}}{k_{1}^{\alpha }+k_{2}^{\alpha }}\left(
k_{1}^{\alpha }N_{1}^{\alpha }+k_{2}^{\alpha }N_{2}^{\alpha }\right) .
\end{equation*}%
Then, the first curvature and the principal normal vector field\ of curve $%
\beta $ are respectively,%
\begin{equation*}
\kappa _{\beta }=\left\Vert T_{\beta }^{\prime }\right\Vert =\frac{\sqrt{2}%
\sqrt{\left( k_{1}^{\alpha }\right) ^{2}+\left( k_{2}^{\alpha }\right) ^{2}}%
}{k_{1}^{\alpha }+k_{2}^{\alpha }}
\end{equation*}%
and%
\begin{equation*}
N_{\beta }=\frac{-1}{\sqrt{\left( k_{1}^{\alpha }\right) ^{2}+\left(
k_{2}^{\alpha }\right) ^{2}}}\left( k_{1}^{\alpha }N_{1}^{\alpha
}+k_{2}^{\alpha }N_{2}^{\alpha }\right) .
\end{equation*}%
On the other hand, we express%
\begin{equation*}
B_{\beta }=\frac{1}{\sqrt{\left( k_{1}^{\alpha }\right) ^{2}+\left(
k_{2}^{\alpha }\right) ^{2}}}\left\vert 
\begin{array}{ccc}
T_{\alpha } & N_{1}^{\alpha } & N_{2}^{\alpha } \\ 
-1 & 0 & 0 \\ 
0 & -k_{1}^{\alpha } & -k_{2}^{\alpha }%
\end{array}%
\right\vert .
\end{equation*}%
So, the binormal vector of curve $\beta $ is%
\begin{equation*}
B_{\beta }=\frac{-1}{\sqrt{\left( k_{1}^{\alpha }\right) ^{2}+\left(
k_{2}^{\alpha }\right) ^{2}}}\left( k_{2}^{\alpha }N_{1}^{\alpha
}-k_{1}^{\alpha }N_{2}^{\alpha }\right) .
\end{equation*}%
We differentiate (\ref{c5}) with respect to $s$ in order to calculate the
torsion of curve $\beta $%
\begin{equation*}
\overset{\cdot \cdot }{\beta }=\frac{-1}{\sqrt{2}}\left\{ \left( \overset{%
\cdot }{k_{1}^{\alpha }}+\overset{\cdot }{k_{2}^{\alpha }}\right) T_{\alpha
}+\left( \left( k_{1}^{\alpha }\right) ^{2}+k_{1}^{\alpha }k_{2}^{\alpha
}\right) N_{1}^{\alpha }+\left( k_{1}^{\alpha }k_{2}^{\alpha }+\left(
k_{2}^{\alpha }\right) ^{2}\right) N_{2}^{\alpha }\right\} ,
\end{equation*}%
and similarly%
\begin{equation*}
\overset{\cdot \cdot \cdot }{\beta }=\frac{-1}{\sqrt{2}}\left( \rho
_{1}T_{\alpha }+\rho _{2}N_{1}^{\alpha }+\rho _{3}N_{2}^{\alpha }\right) ,
\end{equation*}%
where 
\begin{eqnarray*}
\rho _{1} &=&\overset{\cdot \cdot }{k_{1}^{\alpha }}+\overset{\cdot \cdot }{%
k_{2}^{\alpha }}-\left( k_{1}^{\alpha }\right) ^{3}-\left( k_{1}^{\alpha
}\right) ^{2}k_{2}^{\alpha }-k_{1}^{\alpha }\left( k_{2}^{\alpha }\right)
^{2}-\left( k_{2}^{\alpha }\right) ^{3} \\
\rho _{2} &=&3k_{1}^{\alpha }\overset{\cdot }{k_{1}^{\alpha }}%
+2k_{1}^{\alpha }\overset{\cdot }{k_{2}^{\alpha }}+\overset{\cdot }{%
k_{1}^{\alpha }}k_{2}^{\alpha } \\
\rho _{3} &=&2\overset{\cdot }{k_{1}^{\alpha }}k_{2}^{\alpha
}+3k_{2}^{\alpha }\overset{\cdot }{k_{2}^{\alpha }}+k_{1}^{\alpha }\overset{%
\cdot }{k_{2}^{\alpha }}.
\end{eqnarray*}%
The torsion of curve\ $\beta $ is

\begin{equation*}
\tau _{\beta }=\frac{-\sqrt{2}\left \{ \rho _{3}(\left( k_{1}^{\alpha
}\right) ^{2}+k_{1}^{\alpha }k_{2}^{\alpha })-\rho _{2}(k_{1}^{\alpha
}k_{2}^{\alpha }+\left( k_{2}^{\alpha }\right) ^{2})\right \} }{\left(
k_{1}^{\alpha }+k_{2}^{\alpha }\right) \left \{ (k_{1}^{\alpha
}k_{2}^{\alpha }+\left( k_{2}^{\alpha }\right) ^{2})^{2}+(\left(
k_{1}^{\alpha }\right) ^{2}+k_{1}^{\alpha }k_{2}^{\alpha })^{2}\right \} }.
\end{equation*}

The first and second normal field of curve $\beta $ are as follow. Then,
from (\ref{x4}) we obtain

\begin{equation*}
N_{1}^{\beta }=\frac{1}{\sqrt{\left( k_{1}^{\alpha }\right) ^{2}+\left(
k_{2}^{\alpha }\right) ^{2}}}\left \{ \left( k_{2}^{\alpha }\sin \theta
_{\beta }-k_{1}^{\alpha }\cos \theta _{\beta }\right) N_{1}^{\alpha }-\left(
k_{2}^{\alpha }\cos \theta _{\beta }+k_{1}^{\alpha }\sin \theta _{\beta
}\right) N_{2}^{\alpha }\right \} ,
\end{equation*}%
and

\begin{equation*}
N_{2}^{\beta }=\frac{1}{\sqrt{\left( k_{1}^{\alpha }\right) ^{2}+\left(
k_{2}^{\alpha }\right) ^{2}}}\left \{ -\left( k_{1}^{\alpha }\sin \theta
_{\beta }+k_{2}^{\alpha }\cos \theta _{\beta }\right) N_{1}^{\alpha }+\left(
k_{1}^{\alpha }\cos \theta _{\beta }-k_{2}^{\alpha }\sin \theta _{\beta
}\right) N_{2}^{\alpha }\right \}
\end{equation*}%
where $\theta _{\beta }=-\int \tau _{\beta }(s)ds.$

Now, we can calculate\ the natural curvatures of curve $\beta $, so from (%
\ref{x3}) we get

\begin{equation*}
k_{1}^{\beta }=\frac{\sqrt{2(\left( k_{1}^{\alpha }\right) ^{2}+\left(
k_{2}^{\alpha }\right) ^{2})}\cos \theta _{\beta }}{k_{1}^{\alpha
}+k_{2}^{\alpha }}
\end{equation*}%
and

\begin{equation*}
k_{2}^{\beta }=\frac{\sqrt{2(\left( k_{1}^{\alpha }\right) ^{2}+\left(
k_{2}^{\alpha }\right) ^{2})}\sin \theta _{\beta }}{k_{1}^{\alpha
}+k_{2}^{\alpha }}.
\end{equation*}

\subsection{$TN_{1}N_{2}-$Smarandache Curves}

\begin{definition}
Let $\alpha =\alpha (s)$ be a unit speed regular curve in $E^{3}$ and $%
\left\{ T_{\alpha },N_{1}^{\alpha },N_{2}^{\alpha }\right\} $ be its moving
Bishop frame. $T_{\alpha }N_{1}^{\alpha }N_{2}^{\alpha }-$Smarandache curves
can be defined by%
\begin{equation}
\beta (s^{\ast })=\frac{1}{\sqrt{3}}\left( T_{\alpha }+N_{1}^{\alpha
}+N_{2}^{\alpha }\right) .  \label{d1}
\end{equation}
\end{definition}

Now, we can investigate Bishop invariants of $T_{\alpha }N_{1}^{\alpha
}N_{2}^{\alpha }-$Smarandache curves according to $\alpha =\alpha (s).$
Differentiating (\ref{d1}) with respect to $s$, we get%
\begin{equation}
\overset{\cdot }{\beta }=\frac{d\beta }{ds^{\ast }}\frac{ds^{\ast }}{ds}=%
\frac{-1}{\sqrt{3}}\left( \left( k_{1}^{\alpha }+k_{2}^{\alpha }\right)
T_{\alpha }-k_{1}^{\alpha }N_{1}^{\alpha }-k_{2}^{\alpha }N_{2}^{\alpha
}\right)  \label{d5}
\end{equation}%
and%
\begin{equation*}
T_{\beta }\frac{ds^{\ast }}{ds}=\frac{-1}{\sqrt{3}}\left( \left(
k_{1}^{\alpha }+k_{2}^{\alpha }\right) T_{\alpha }-k_{1}^{\alpha
}N_{1}^{\alpha }-k_{2}^{\alpha }N_{2}^{\alpha }\right)
\end{equation*}%
where%
\begin{equation}
\frac{ds^{\ast }}{ds}=\sqrt{\frac{2\left( \left( k_{1}^{\alpha }\right)
^{2}+k_{1}^{\alpha }k_{2}^{\alpha }+\left( k_{2}^{\alpha }\right)
^{2}\right) }{3}}.  \label{d3}
\end{equation}%
The tangent vector of curve $\beta $ can be written as follow,%
\begin{equation}
T_{\beta }=\frac{-1}{\sqrt{2\left\{ \left( k_{1}^{\alpha }\right)
^{2}+k_{1}^{\alpha }k_{2}^{\alpha }+\left( k_{2}^{\alpha }\right)
^{2}\right\} }}\left( \left( k_{1}^{\alpha }+k_{2}^{\alpha }\right)
T_{\alpha }-k_{1}^{\alpha }N_{1}^{\alpha }-k_{2}^{\alpha }N_{2}^{\alpha
}\right) .  \label{d2}
\end{equation}%
Differentiating (\ref{d2}) with respect to $s$, we obtain%
\begin{equation}
\frac{dT_{\beta }}{ds^{\ast }}\frac{ds^{\ast }}{ds}=\frac{1}{2\sqrt{2}%
\left\{ \left( k_{1}^{\alpha }\right) ^{2}+k_{1}^{\alpha }k_{2}^{\alpha
}+\left( k_{2}^{\alpha }\right) ^{2}\right\} ^{\frac{3}{2}}}\left( \lambda
_{1}T_{\alpha }+\lambda _{2}N_{1}^{\alpha }+\lambda _{3}N_{2}^{\alpha
}\right)  \label{d4}
\end{equation}%
where%
\begin{eqnarray*}
\lambda _{1} &=&-2\overset{\cdot }{k_{1}^{\alpha }}\left( k_{2}^{\alpha
}\right) ^{2}-\left( k_{1}^{\alpha }\right) ^{2}\overset{\cdot }{%
k_{2}^{\alpha }}+k_{1}^{\alpha }k_{2}^{\alpha }\overset{\cdot }{%
k_{2}^{\alpha }}-2\left( k_{1}^{\alpha }\right) ^{4}-2\left( k_{1}^{\alpha
}\right) ^{3}k_{2}^{\alpha } \\
&&-4\left( k_{1}^{\alpha }\right) ^{2}\left( k_{2}^{\alpha }\right)
^{2}-2k_{1}^{\alpha }\left( k_{2}^{\alpha }\right) ^{3}-2\left(
k_{2}^{\alpha }\right) ^{4}+k_{1}^{\alpha }\overset{\cdot }{k_{1}^{\alpha }}%
k_{2}^{\alpha }+\overset{\cdot }{k_{1}^{\alpha }}\left( k_{2}^{\alpha
}\right) ^{2} \\
\lambda _{2} &=&-2\left( k_{1}^{\alpha }\right) ^{4}-4\left( k_{1}^{\alpha
}\right) ^{3}k_{2}^{\alpha }-4\left( k_{1}^{\alpha }\right) ^{2}\left(
k_{2}^{\alpha }\right) ^{2}-2k_{1}^{\alpha }\left( k_{2}^{\alpha }\right)
^{3}+k_{1}^{\alpha }\overset{\cdot }{k_{1}^{\alpha }}k_{2}^{\alpha } \\
&&+2\overset{\cdot }{k_{1}^{\alpha }}\left( k_{2}^{\alpha }\right)
^{2}-\left( k_{1}^{\alpha }\right) ^{2}\overset{\cdot }{k_{2}^{\alpha }}%
-2k_{1}^{\alpha }k_{2}^{\alpha }\overset{\cdot }{k_{2}^{\alpha }} \\
\lambda _{3} &=&-2\left( k_{1}^{\alpha }\right) ^{3}k_{2}^{\alpha }-4\left(
k_{1}^{\alpha }\right) ^{2}\left( k_{2}^{\alpha }\right) ^{2}-4k_{1}^{\alpha
}\left( k_{2}^{\alpha }\right) ^{3}-2\left( k_{2}^{\alpha }\right) ^{4} \\
&&+2\left( k_{1}^{\alpha }\right) ^{2}\overset{\cdot }{k_{2}^{\alpha }}%
+k_{1}^{\alpha }k_{2}^{\alpha }\overset{\cdot }{k_{2}^{\alpha }}%
-2k_{1}^{\alpha }\overset{\cdot }{k_{1}^{\alpha }}k_{2}^{\alpha }-\overset{%
\cdot }{k_{1}^{\alpha }}\left( k_{2}^{\alpha }\right) ^{2}.
\end{eqnarray*}%
Substituting (\ref{d3}) in (\ref{d4}), we get%
\begin{equation*}
T_{\beta }^{\prime }=\frac{\sqrt{3}}{4\left\{ \left( k_{1}^{\alpha }\right)
^{2}+k_{1}^{\alpha }k_{2}^{\alpha }+\left( k_{2}^{\alpha }\right)
^{2}\right\} ^{2}}\left( \lambda _{1}T_{\alpha }+\lambda _{2}N_{1}^{\alpha
}+\lambda _{3}N_{2}^{\alpha }\right) .
\end{equation*}%
Then, the first curvature and the principal normal vector field of curve $%
\beta $ are respectively,%
\begin{equation*}
\kappa _{\beta }=\left\Vert T_{\beta }^{\prime }\right\Vert =\frac{\sqrt{3}%
\sqrt{\lambda _{1}^{2}+\lambda _{2}^{2}+\lambda _{3}^{2}}}{4\left\{ \left(
k_{1}^{\alpha }\right) ^{2}+k_{1}^{\alpha }k_{2}^{\alpha }+\left(
k_{2}^{\alpha }\right) ^{2}\right\} ^{2}}
\end{equation*}%
and%
\begin{equation}
N_{\beta }=\frac{1}{\sqrt{\lambda _{1}^{2}+\lambda _{2}^{2}+\lambda _{3}^{2}}%
}\left( \lambda _{1}T_{\alpha }+\lambda _{2}N_{1}^{\alpha }+\lambda
_{3}N_{2}^{\alpha }\right) .  \label{d6}
\end{equation}%
On the other hand, we express%
\begin{equation*}
B_{\beta }=\frac{1}{\sqrt{2\left( \left( k_{1}^{\alpha }\right)
^{2}+k_{1}^{\alpha }k_{2}^{\alpha }+\left( k_{2}^{\alpha }\right)
^{2}\right) }\sqrt{\lambda _{1}^{2}+\lambda _{2}^{2}+\lambda _{3}^{2}}}%
\left\vert 
\begin{array}{ccc}
T_{\alpha } & N_{1}^{\alpha } & N_{2}^{\alpha } \\ 
-\left( k_{1}^{\alpha }+k_{2}^{\alpha }\right) & k_{1}^{\alpha } & 
k_{2}^{\alpha } \\ 
\lambda _{1} & \lambda _{2} & \lambda _{3}%
\end{array}%
\right\vert .
\end{equation*}%
So, the binormal vector of curve $\beta $ is%
\begin{equation}
B_{\beta }=\frac{1}{\sqrt{2\left( \left( k_{1}^{\alpha }\right)
^{2}+k_{1}^{\alpha }k_{2}^{\alpha }+\left( k_{2}^{\alpha }\right)
^{2}\right) }\sqrt{\lambda _{1}^{2}+\lambda _{2}^{2}+\lambda _{3}^{2}}}%
\left( \sigma _{1}T_{\alpha }+\sigma _{2}N_{1}^{\alpha }+\sigma
_{3}N_{2}^{\alpha }\right)  \label{d7}
\end{equation}%
where%
\begin{equation*}
\sigma _{1}=\lambda _{3}k_{1}^{\alpha }-\lambda _{2}k_{2}^{\alpha }\text{ \
, \ }\sigma _{2}=\lambda _{3}\left( k_{1}^{\alpha }+k_{2}^{\alpha }\right)
+\lambda _{1}k_{2}^{\alpha }\text{ \ , \ }\sigma _{3}=-\left( \lambda
_{2}\left( k_{1}^{\alpha }+k_{2}^{\alpha }\right) +\lambda _{1}k_{1}^{\alpha
}\right) .
\end{equation*}%
We differentiate (\ref{d5}) with respect to $s$ in order to calculate the
torsion of curve $\beta ,$%
\begin{equation*}
\overset{\cdot \cdot }{\beta }=\frac{-1}{\sqrt{3}}\left\{ (\overset{\cdot }{%
k_{1}^{\alpha }}+\overset{\cdot }{k_{2}^{\alpha }}+\left( k_{1}^{\alpha
}\right) ^{2}+\left( k_{2}^{\alpha }\right) ^{2})T_{\alpha }-(\overset{\cdot 
}{k_{1}^{\alpha }}-\left( k_{1}^{\alpha }\right) ^{2}-k_{1}^{\alpha
}k_{2}^{\alpha })N_{1}^{\alpha }-(\overset{\cdot }{k_{2}^{\alpha }}%
-k_{1}^{\alpha }k_{2}^{\alpha }-\left( k_{2}^{\alpha }\right)
^{2})N_{2}^{\alpha }\right\} ,
\end{equation*}%
and similarly%
\begin{equation*}
\overset{\cdot \cdot \cdot }{\beta }=\frac{1}{\sqrt{3}}\left( \rho
_{1}T_{\alpha }+\rho _{2}N_{1}^{\alpha }+\rho _{3}N_{2}^{\alpha }\right) ,
\end{equation*}%
where 
\begin{eqnarray*}
\rho _{1} &=&-\overset{\cdot \cdot }{k_{1}^{\alpha }}-\overset{\cdot \cdot }{%
k_{2}^{\alpha }}-3k_{1}^{\alpha }\overset{\cdot }{k_{1}^{\alpha }}%
-3k_{2}^{\alpha }\overset{\cdot }{k_{2}^{\alpha }}+\left( k_{1}^{\alpha
}\right) ^{3}+\left( k_{1}^{\alpha }\right) ^{2}k_{2}^{\alpha
}+k_{1}^{\alpha }\left( k_{2}^{\alpha }\right) ^{2}+\left( k_{2}^{\alpha
}\right) ^{3} \\
\rho _{2} &=&-3k_{1}^{\alpha }\overset{\cdot }{k_{1}^{\alpha }}%
-2k_{1}^{\alpha }\overset{\cdot }{k_{2}^{\alpha }}-\left( k_{1}^{\alpha
}\right) ^{3}-k_{1}^{\alpha }\left( k_{2}^{\alpha }\right) ^{2}-\overset{%
\cdot }{k_{1}^{\alpha }}k_{2}^{\alpha }+\overset{\cdot \cdot }{k_{1}^{\alpha
}} \\
\rho _{3} &=&-2\overset{\cdot }{k_{1}^{\alpha }}k_{2}^{\alpha
}-3k_{2}^{\alpha }\overset{\cdot }{k_{2}^{\alpha }}-\left( k_{1}^{\alpha
}\right) ^{2}k_{2}^{\alpha }-\left( k_{2}^{\alpha }\right)
^{3}-k_{1}^{\alpha }\overset{\cdot }{k_{2}^{\alpha }}+\overset{\cdot \cdot }{%
k_{2}^{\alpha }}
\end{eqnarray*}%
The torsion of curve\ $\beta $ is

\begin{equation*}
\tau _{\beta }=\frac{\sqrt{3}\left \{ 
\begin{array}{c}
(\overset{\cdot }{k_{1}^{\alpha }}-\left( k_{1}^{\alpha }\right)
^{2}-k_{1}^{\alpha }k_{2}^{\alpha })(-\rho _{3}k_{1}^{\alpha }-\rho
_{3}k_{2}^{\alpha }-\rho _{1}k_{2}^{\alpha }) \\ 
+(\overset{\cdot }{k_{2}^{\alpha }}-k_{1}^{\alpha }k_{2}^{\alpha }-\left(
k_{2}^{\alpha }\right) ^{2})(\rho _{2}k_{1}^{\alpha }+\rho _{2}k_{2}^{\alpha
}+\rho _{1}k_{1}^{\alpha }) \\ 
-(\overset{\cdot }{k_{1}^{\alpha }}+\overset{\cdot }{k_{2}^{\alpha }}+\left(
k_{1}^{\alpha }\right) ^{2}+\left( k_{2}^{\alpha }\right) ^{2})(-\rho
_{3}k_{1}^{\alpha }+\rho _{2}k_{2}^{\alpha })%
\end{array}%
\right \} }{\left \{ 
\begin{array}{c}
(k_{1}^{\alpha }\overset{\cdot }{k_{2}^{\alpha }}-\overset{\cdot }{%
k_{1}^{\alpha }}k_{2}^{\alpha })^{2}+(k_{1}^{\alpha }\overset{\cdot }{%
k_{2}^{\alpha }}-2\left( k_{1}^{\alpha }\right) ^{2}k_{2}^{\alpha
}-2k_{1}^{\alpha }\left( k_{2}^{\alpha }\right) ^{2}-2\left( k_{2}^{\alpha
}\right) ^{3}-\overset{\cdot }{k_{1}^{\alpha }}k_{2}^{\alpha })^{2} \\ 
+(2\left( k_{1}^{\alpha }\right) ^{3}+2\left( k_{1}^{\alpha }\right)
^{2}k_{2}^{\alpha }+2k_{1}^{\alpha }\left( k_{2}^{\alpha }\right) ^{2}-%
\overset{\cdot }{k_{1}^{\alpha }}k_{2}^{\alpha }+k_{1}^{\alpha }\overset{%
\cdot }{k_{2}^{\alpha }})^{2}%
\end{array}%
\right \} }.
\end{equation*}

The first and second normal field of curve $\beta $ are as follow. Then,
from (\ref{x4}) we obtain

\begin{equation*}
N_{1}^{\beta }=\frac{1}{\sqrt{2\mu }\eta }\left \{ 
\begin{array}{c}
\left( \sqrt{2\mu }\cos \theta _{\beta }\lambda _{1}-\sin \theta _{\beta
}\sigma _{1}\right) T_{\alpha } \\ 
+\left( \sqrt{2\mu }\cos \theta _{\beta }\lambda _{2}-\sin \theta _{\beta
}\sigma _{2}\right) N_{1}^{\alpha } \\ 
+\left( \sqrt{2\mu }\cos \theta _{\beta }\lambda _{3}-\sin \theta _{\beta
}\sigma _{3}\right) N_{2}^{\alpha }%
\end{array}%
\right \} ,
\end{equation*}%
and

\begin{equation*}
N_{2}^{\beta }=\frac{1}{\sqrt{2\mu }\eta }\left \{ 
\begin{array}{c}
\left( \sqrt{2\mu }\sin \theta _{\beta }\lambda _{1}+\cos \theta _{\beta
}\sigma _{1}\right) T_{\alpha } \\ 
+\left( \sqrt{2\mu }\sin \theta _{\beta }\lambda _{2}+\cos \theta _{\beta
}\sigma _{2}\right) N_{1}^{\alpha } \\ 
+\left( \sqrt{2\mu }\sin \theta _{\beta }\lambda _{3}+\cos \theta _{\beta
}\sigma _{3}\right) N_{2}^{\alpha }%
\end{array}%
\right \} ,
\end{equation*}%
where $\theta _{\beta }=-\int \tau _{\beta }(s)ds$ and $\mu =\left(
k_{1}^{\alpha }\right) ^{2}+k_{1}^{\alpha }k_{2}^{\alpha }+\left(
k_{2}^{\alpha }\right) ^{2}$, $\eta =\sqrt{\lambda _{1}^{2}+\lambda
_{2}^{2}+\lambda _{3}^{2}}$

Now, we can calculate\ the natural curvatures of curve $\beta $, so from (%
\ref{x3}) we get

\begin{equation*}
k_{1}^{\beta }=\frac{\sqrt{3(\lambda _{1}^{2}+\lambda _{2}^{2}+\lambda
_{3}^{2})}\cos \theta _{\beta }}{4(\left( k_{1}^{\alpha }\right)
^{2}+k_{1}^{\alpha }k_{2}^{\alpha }+\left( k_{2}^{\alpha }\right) ^{2})^{2}}
\end{equation*}%
and

\begin{equation*}
k_{2}^{\beta }=\frac{\sqrt{3(\lambda _{1}^{2}+\lambda _{2}^{2}+\lambda
_{3}^{2})}\sin \theta _{\beta }}{4(\left( k_{1}^{\alpha }\right)
^{2}+k_{1}^{\alpha }k_{2}^{\alpha }+\left( k_{2}^{\alpha }\right) ^{2})^{2}}.
\end{equation*}

\section{Curvature Spheres and Osculating Spheres of Smarandache Curves}

\begin{definition}
Let $\alpha :(a,b)\rightarrow 
\mathbb{R}
^{n}$ be a regular curve and let $F:%
\mathbb{R}
^{n}\rightarrow 
\mathbb{R}
$ be a differentiable function. We say that a parametrically defined curve $%
\alpha $ and an implicity defined hypersurface $F^{-1}(0)$ have contact of
order $k$ at $\alpha (t_{0})$ provided we have $(Fo\alpha )(t_{0})=(Fo\alpha
)^{\prime }(t_{0})=...=(Fo\alpha )^{(k)}(t_{0})=0$ but $(Fo\alpha
)^{(k+1)}(t_{0})\neq 0$ [5].
\end{definition}

If any circle which has a contact of first-order with a planar curve, it is
called the \textit{curvature circle}. If any sphere which has a contact of
second-order with a space curve, It is called the \textit{curvature sphere}.
Furthermore, if any sphere which has a contact of third-order with a space
curve, it is called the \textit{osculating sphere}.

\begin{theorem}
Let $\beta (s^{\ast })$ be a unit speed Smarandache curve of $\alpha (s)$
with first curvature $k_{1}^{\beta }$ and second curvature $k_{2}^{\beta }$,
then the centres of curvature spheres of $\beta (s^{\ast })$ according to
Bishop frame are%
\begin{eqnarray*}
c(s^{\ast }) &=&\beta (s^{\ast })+\frac{k_{1}^{\beta }\mp k_{2}^{\beta }%
\sqrt{r^{2}\left( (k_{1}^{\beta })^{2}-(k_{2}^{\beta })^{2}\right) -1}}{%
4(k_{1}^{\beta })^{2}}N_{1}^{\beta } \\
&&+\frac{3k_{1}^{\beta }\pm k_{2}^{\beta }\sqrt{r^{2}((k_{1}^{\beta
})^{2}-(k_{2}^{\beta })^{2})-1}}{4k_{1}^{\beta }k_{2}^{\beta }}N_{2}^{\beta
}.
\end{eqnarray*}
\end{theorem}

\begin{proof}
Assume that $c(s^{\ast })$ is the centre of curvature sphere of $\beta
(s^{\ast })$, then radius vector of curvature sphere is 
\begin{equation*}
c-\beta =\delta _{1}T_{\beta }+\delta _{2}N_{1}^{\beta }+\delta
_{3}N_{2}^{\beta }.
\end{equation*}

Let we define the function $F$ such that $F=\left\langle c(s^{\ast })-\beta
(s^{\ast }),c(s^{\ast })-\beta (s^{\ast })\right\rangle -r^{2}$. First and
second derivatives of $F$ are as follow.%
\begin{equation*}
F^{\prime }=-2\delta _{1}
\end{equation*}%
\begin{equation*}
F^{\prime \prime }=-2(k_{1}^{\beta }\delta _{2}+k_{2}^{\beta }\delta _{3}-1).
\end{equation*}%
The condition of contact of second-order requires $F=0,$ $F^{\prime }=0$ and 
$F^{\prime \prime }=0$. Then, coefficients $\delta _{1},\delta _{2},\delta
_{3}$ are obtained%
\begin{eqnarray*}
\delta _{1} &=&0 \\
\delta _{2} &=&\frac{k_{1}^{\beta }\mp k_{2}^{\beta }\sqrt{r^{2}\left(
(k_{1}^{\beta })^{2}-(k_{2}^{\beta })^{2}\right) -1}}{4(k_{1}^{\beta })^{2}}
\\
\delta _{3} &=&\frac{3k_{1}^{\beta }\pm k_{2}^{\beta }\sqrt{%
r^{2}((k_{1}^{\beta })^{2}-(k_{2}^{\beta })^{2})-1}}{4k_{1}^{\beta
}k_{2}^{\beta }}.
\end{eqnarray*}%
Thus, centres of the curvature spheres of Smarandache curve $\beta $ are%
\begin{eqnarray*}
c(s^{\ast }) &=&\beta (s^{\ast })+\frac{k_{1}^{\beta }\mp k_{2}^{\beta }%
\sqrt{r^{2}\left( (k_{1}^{\beta })^{2}-(k_{2}^{\beta })^{2}\right) -1}}{%
4(k_{1}^{\beta })^{2}}N_{1}^{\beta } \\
&&+\frac{3k_{1}^{\beta }\pm k_{2}^{\beta }\sqrt{r^{2}((k_{1}^{\beta
})^{2}-(k_{2}^{\beta })^{2})-1}}{4k_{1}^{\beta }k_{2}^{\beta }}N_{2}^{\beta
}.
\end{eqnarray*}
\end{proof}

Let $\delta _{3}=\lambda ,$ then $c(s^{\ast })=\beta (s^{\ast })+\delta
_{2}N_{1}^{\beta }+\lambda N_{2}^{\beta }$ is a line passing the point $%
\beta (s^{\ast })+\delta _{2}N_{1}^{\beta }$ and parallel to line $%
N_{2}^{\beta }.$ Thus we have the following corollary.

\begin{corollary}
Let $\beta (s^{\ast })$ be a unit speed Smarandache curve of $\alpha (s)$
then each centres of the curvature spheres of $\beta $ are on a line.
\end{corollary}

\begin{theorem}
Let $\beta (s^{\ast })$ be a unit speed Smarandache curve of $\alpha (s)$
with first curvature $k_{1}^{\beta }$ and second curvature $k_{2}^{\beta }$
, then the centres of osculating spheres of $\beta (s^{\ast })$ according to
Bishop frame are%
\begin{equation*}
c(s^{\ast })=\beta (s^{\ast })+\frac{(k_{2}^{\beta })^{\prime }}{%
(k_{1}^{\beta })^{2}\left( \frac{k_{2}^{\beta }}{k_{1}^{\beta }}\right)
^{\prime }}N_{1}^{\beta }-\frac{(k_{1}^{\beta })^{\prime }}{(k_{1}^{\beta
})^{2}\left( \frac{k_{2}^{\beta }}{k_{1}^{\beta }}\right) ^{\prime }}%
N_{2}^{\beta }.
\end{equation*}
\end{theorem}

\begin{proof}
Assume that $c(s^{\ast })$ is the centre of osculating sphere of $\beta
(s^{\ast })$, then radius vector of osculating sphere is 
\begin{equation*}
c(s^{\ast })-\beta (s^{\ast })=\delta _{1}T_{\beta }+\delta _{2}N_{1}^{\beta
}+\delta _{3}N_{2}^{\beta }.
\end{equation*}

Let we define the function $F$ such that $F=\left\langle c(s^{\ast })-\beta
(s^{\ast }),c(s^{\ast })-\beta (s^{\ast })\right\rangle -r^{2}$. First,
second and third derivatives of $F$ are as follow.%
\begin{equation*}
F^{\prime }=-2\delta _{1}
\end{equation*}%
\begin{equation*}
F^{\prime \prime }=-2(k_{1}^{\beta }\delta _{2}+k_{2}^{\beta }\delta _{3}-1)
\end{equation*}%
\begin{equation*}
F^{\prime \prime \prime }=-2\left( (k_{1}^{\beta })^{\prime }\delta
_{2}-(k_{1}^{\beta })^{2}\delta _{1}+(k_{2}^{\beta })^{\prime }\delta
_{3}-(k_{2}^{\beta })^{2}\delta _{1}\right) .
\end{equation*}%
The condition of contact of third-order requires $F=0,$ $F^{\prime }=0,$ $%
F^{\prime \prime }=0$ and $F^{\prime \prime \prime }=0$. Then, coefficients $%
\delta _{1},\delta _{2},\delta _{3}$ are obtained%
\begin{equation*}
\delta _{1}=0,\text{ \ \ \ }\delta _{2}=\frac{(k_{2}^{\beta })^{\prime }}{%
(k_{1}^{\beta })^{2}\left( \frac{k_{2}^{\beta }}{k_{1}^{\beta }}\right)
^{\prime }},\text{ \ \ }\delta _{3}=-\frac{(k_{1}^{\beta })^{\prime }}{%
(k_{1}^{\beta })^{2}\left( \frac{k_{2}^{\beta }}{k_{1}^{\beta }}\right)
^{\prime }}.
\end{equation*}%
Thus, centres of the osculating spheres of Smarandache curve $\beta $ can be
written as follow%
\begin{equation*}
c(s^{\ast })=\beta (s^{\ast })+\frac{(k_{2}^{\beta })^{\prime }}{%
(k_{1}^{\beta })^{2}\left( \frac{k_{2}^{\beta }}{k_{1}^{\beta }}\right)
^{\prime }}N_{1}^{\beta }-\frac{(k_{1}^{\beta })^{\prime }}{(k_{1}^{\beta
})^{2}\left( \frac{k_{2}^{\beta }}{k_{1}^{\beta }}\right) ^{\prime }}%
N_{2}^{\beta }.
\end{equation*}
\end{proof}

\begin{corollary}
\bigskip Let $\beta (s^{\ast })$ be a unit speed Smarandache curve of $%
\alpha (s)$ with first curvature $k_{1}^{\beta }$ and second curvature $%
k_{2}^{\beta }$ , then the radius of the osculating sphere of $\beta
(s^{\ast })$ is%
\begin{equation*}
r(s^{\ast })=\frac{\sqrt{\left( (k_{1}^{\beta })^{\prime }\right)
^{2}+\left( (k_{2}^{\beta })^{\prime }\right) ^{2}}}{(k_{1}^{\beta
})^{2}\left( \frac{k_{2}^{\beta }}{k_{1}^{\beta }}\right) ^{\prime }}.
\end{equation*}%
\bigskip
\end{corollary}

\begin{example}
Let $\alpha \left( t\right) $ be the Salkowski curve given by%
\begin{equation}
\alpha \left( t\right) =\left( \alpha _{1}\left( t\right) ,\alpha _{2}\left(
t\right) ,\alpha _{3}\left( t\right) \right)   \label{x2}
\end{equation}%
where%
\begin{eqnarray*}
\alpha _{1} &=&\frac{1}{\sqrt{1+m^{2}}}\left( -\frac{1-n}{4\left(
1+2n\right) }\sin \left( \left( 1+2n\right) t\right) -\frac{1+n}{4\left(
1-2n\right) }\sin \left( \left( 1-2n\right) t\right) -\frac{1}{2}\sin
t\right)  \\
\alpha _{2} &=&\frac{1}{\sqrt{1+m^{2}}}\left( \frac{1-n}{4\left( 1+2n\right) 
}\cos \left( \left( 1+2n\right) t\right) +\frac{1+n}{4\left( 1-2n\right) }%
\cos \left( \left( 1-2n\right) t\right) +\frac{1}{2}\cos t\right)  \\
\alpha _{3} &=&\frac{\cos \left( 2nt\right) }{4m\sqrt{1+m^{2}}}
\end{eqnarray*}%
and%
\begin{equation*}
n=\frac{m}{\sqrt{1+m^{2}}}.
\end{equation*}%
$\alpha \left( t\right) $ can be expressed with the arc lenght parameter. We
get unit speed regular curve by using the parametrization $t=\frac{1}{n}%
\arcsin \left( n\sqrt{1+m^{2}}s\right) $. Furthermore, graphics of special
Smarandache curves are as follow. Here $m=\sqrt{3}.$%
\begin{equation*}
\FRAME{itbpF}{1.9925in}{2.1586in}{0in}{}{}{fig11.png}{\special{language
"Scientific Word";type "GRAPHIC";maintain-aspect-ratio TRUE;display
"USEDEF";valid_file "F";width 1.9925in;height 2.1586in;depth
0in;original-width 3.0831in;original-height 3.3434in;cropleft "0";croptop
"1";cropright "1";cropbottom "0";filename 'fig11.png';file-properties
"XNPEU";}}
\end{equation*}%
\begin{equation*}
\text{Figure 1 : Smarandache Curve }T_{\alpha }N_{1}^{\alpha }
\end{equation*}%
\begin{equation*}
\FRAME{itbpF}{2.0159in}{1.8481in}{0in}{}{}{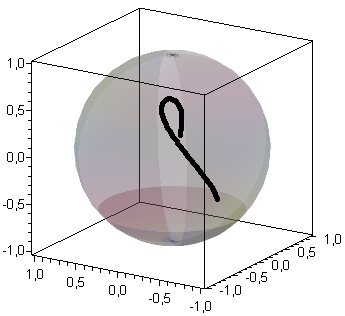}{\special{language
"Scientific Word";type "GRAPHIC";maintain-aspect-ratio TRUE;display
"USEDEF";valid_file "F";width 2.0159in;height 1.8481in;depth
0in;original-width 3.5942in;original-height 3.2915in;cropleft "0";croptop
"1";cropright "1";cropbottom "0";filename 'fig22.png';file-properties
"XNPEU";}}
\end{equation*}%
\begin{equation*}
\text{Figure 2 : Smarandache Curve }T_{\alpha }N_{2}^{\alpha }
\end{equation*}%
\begin{equation*}
\FRAME{itbpF}{2.1404in}{1.8317in}{0in}{}{}{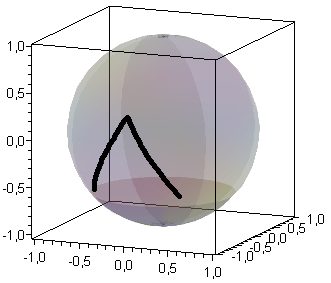}{\special{language
"Scientific Word";type "GRAPHIC";maintain-aspect-ratio TRUE;display
"USEDEF";valid_file "F";width 2.1404in;height 1.8317in;depth
0in;original-width 3.448in;original-height 2.9482in;cropleft "0";croptop
"1";cropright "1";cropbottom "0";filename 'fig33.png';file-properties
"XNPEU";}}
\end{equation*}%
\begin{equation*}
\text{Figure 3 : Smarandache Curve }N_{1}^{\alpha }N_{2}^{\alpha }
\end{equation*}%
\begin{equation*}
\FRAME{itbpF}{2.098in}{1.7876in}{0in}{}{}{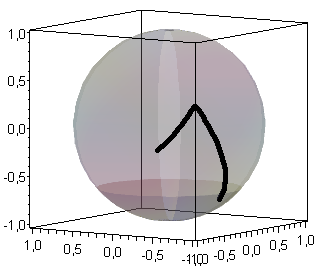}{\special{language
"Scientific Word";type "GRAPHIC";maintain-aspect-ratio TRUE;display
"USEDEF";valid_file "F";width 2.098in;height 1.7876in;depth
0in;original-width 3.3434in;original-height 2.8435in;cropleft "0";croptop
"1";cropright "1";cropbottom "0";filename 'fig44.png';file-properties
"XNPEU";}}
\end{equation*}%
\begin{equation*}
\text{Figure 4 : Smarandache Curve }T_{\alpha }N_{1}^{\alpha }N_{2}^{\alpha }
\end{equation*}
\end{example}


\begin{thebibliography}{99}
\bibitem{1} Ali, \ A. T., "Special Smarandache Curves in the Euclidean
Space", International Journal of Mathematical Combinatorics, Vol.2,
pp.30-36, 2010.

\bibitem{2} Bishop, L.R., \textquotedblleft There is more than one way to
frame a curve\textquotedblright , Amer. Math. Monthly, Vol.82, Issue.3,
pp.246-251, 1975.

\bibitem{3} Bukcu, B. and Karacan, M.K., "On The Slant Helices According to
Bishop Frame", International Journal of Computational and Mathematical
Sciences, Vol 3, Number 2, Spring 2009.

\bibitem{4} Bukcu, B. and Karacan, M.K., "Special Bishop Motion and Bishop
Darboux Rotation Axis of The Space Curve", Journal of Dynamical Systems and
Geometric Theories,Vol.6, 2008.

\bibitem{5} Gray, A., "Modern Differential Geometry of Curves and Surfaces
with Mathematica", 2nd edition, CRC Press, 1998.

\bibitem{6} Hanson, A. J., and Ma, H., \textquotedblleft Parallel Transport
Approach to Curve Framing\textquotedblright , Indiana University,
Techreports-TR425, January 11, 1995.

\bibitem{7} Hanson, A. J., and Ma, H., \textquotedblleft Quaternion Frame
Approach to Streamline Visualization\textquotedblright , Ieee Transactions
On Visualization And Computer Graphics, Vol.I, No.2, June 1995.

\bibitem{8} Turgut, M. and Yilmaz, S., "Smarandache Curves in Minkowski
Space-time", International Journal of Mathematical Combinatorics, Vol.3,
pp.51-55, 2008.

\bibitem{9} Turhan, E. and K\"{o}rp\i nar, T., "Biharmonic Slant Helices
According to Bishop Frame in $E^{3}$", International Journal of Mathematical
Combinatorics, Vol.3, pp.64-68, 2010.

\bibitem{10} Yilmaz, S. and Turgut, M., "On the Differential Geometry of the
Curves in Minkowski Space-time I", Int. J. Contemp. Math. Sciences, Vol.3,
No.27, pp.1343-1349, 2008.
\end{thebibliography}
\end{document}